\documentclass[11pt]{amsart}
\usepackage{amsxtra}
\theoremstyle{plain}

\newcommand{\cleqn}{\setcounter{equation}{0}}
\newcommand{\clth}{\setcounter{theorem}{0}}
\newcommand {\sectionnew}[1]{\section{#1}\cleqn\clth}

\newtheorem{theorem}{Theorem}[section]
\newtheorem{lemma}[theorem]{Lemma}
\newtheorem{definition-theorem}[theorem]{Definition-Theorem}
\newtheorem{proposition}[theorem]{Proposition}
\newtheorem{corollary}[theorem]{Corollary}
\newtheorem{definition}[theorem]{Definition}
\newtheorem{example}[theorem]{Example}
\newtheorem{remark}[theorem]{Remark}
\newtheorem{notation}[theorem]{Notation}
\newcommand \bth[1] { \begin{theorem}\label{t#1} }
\newcommand \ble[1] { \begin{lemma}\label{l#1} }

\newcommand \bpr[1] { \begin{proposition}\label{p#1} }
\newcommand \bco[1] { \begin{corollary}\label{c#1} }
\newcommand \bde[1] { \begin{definition}\label{d#1}\rm }
\newcommand \bex[1] { \begin{example}\label{e#1}\rm }
\newcommand \bre[1] { \begin{remark}\label{r#1}\rm }

\newcommand \bnota[1] { \begin{notation}\label{n#1}\rm }
\newcommand {\eth} { \end{theorem} }
\newcommand {\ele} { \end{lemma} }

\newcommand {\epr} { \end{proposition} }
\newcommand {\eco} { \end{corollary} }
\newcommand {\ede} { \end{definition} }
\newcommand {\eex} { \end{example} }
\newcommand {\ere} { \end{remark} }

\newcommand {\enota} { \end{notation} }
\newcommand \thref[1]{Theorem \ref{t#1}}
\newcommand \leref[1]{Lemma \ref{l#1}}
\newcommand \prref[1]{Proposition \ref{p#1}}

\newcommand \lb[1]{\label{#1}}

\def \Cset {{\mathbb C}}

\def \De {\Delta}   

\def \al {\alpha}

\def \mt  {\mapsto}

\def \hra {\hookrightarrow}

\def \ci  {\circ}

\def \sign { {\mathrm{sign}} }

\def \const { {\mathrm{const}} }

\def \gl  {\mathfrak{gl}}

\def \glnC {{\mathfrak{gl}}_n(\Cset)}
\def \GL  {\mathrm{GL}}
\def \GLn {{\mathrm{GL}}_n}
\def \GLnC {{\mathrm{GL}}_n(\Cset)}

\DeclareMathOperator \Gr { {\mathrm{Gr}} }

\DeclareMathOperator \Vect { {\mathrm{Vect}} }

\begin{document}
\title[Completeness of determinantal Hamiltonian flows]
{Completeness of determinantal Hamiltonian flows on the matrix 
affine Poisson space}
\author[Michael Gekhtman]{Michael Gekhtman}
\address{
Department of Mathematics \\
University of Notre Dame \\
Notre Dame, IN 46556 \\
U.S.A.
}
\email{Michael.Gekhtman.1@nd.edu}
\author[Milen Yakimov]{Milen Yakimov}
\address{ Department of Mathematics \\
Louisiana State Univerity \\
Baton Rouge, LA 70803 and
Department of Mathematics \\
University of California \\
Santa Barbara, CA 93106 \\
U.S.A.
}
\email{yakimov@math.lsu.edu}
\date{}
\begin{abstract} The matrix affine Poisson space $(M_{m,n}, \pi_{m,n})$ 
is the space of complex rectangular matrices equipped with a canonical 
quadratic Poisson structure 
which in the square case $m=n$ reduces to the standard Poisson structure
on $\GL_n(\Cset)$. We prove that the Hamiltonian flows 
of all minors are complete. As a corollary we obtain that all 
Kogan--Zelevinsky integrable systems on $M_{n,n}$ are complete and thus 
induce (analytic) Hamiltonian actions of $\Cset^{n(n-1)/2}$ on 
$(M_{n,n}, \pi_{n,n})$ (as well as on $\GL_n(\Cset)$ 
and on ${\mathrm{SL}}_n(\Cset)$).

We define Gelfand--Zeitlin integrable systems on $(M_{n,n}, \pi_{n,n})$ 
from chains of Poisson projections and prove that their flows are also complete. 
This is an analog for the quadratic Poisson structure $\pi_{n,n}$ 
of the recent result of Kostant and Wallach \cite{KW} that the flows 
of the complexified classical Gelfand--Zeitlin integrable systems are 
complete. 
\end{abstract}
\maketitle
\sectionnew{Introduction}\lb{Introduction}
The matrix affine Poisson space is the complex affine space $M_{m,n}$
consisting of rectangular matrices of size $m \times n$
equipped with the quadratic Poisson structure  
\begin{equation}
\lb{pi_mn}
\pi_{m,n} = \sum_{i, k=1}^m \sum_{j, l=1}^n
(\sign(k-i) + \sign(l-j)) x_{il} x_{kj}
\frac{\partial}{\partial x_{ij}} \wedge
\frac{\partial}{\partial x_{kl}},
\end{equation}
where $x_{ij}$ are the standard coordinate functions on $M_{m,n}.$
The geometry of this Poisson structure was studied
in \cite{GSV1,BGY, GY2}, motivated by its relation to cluster
algebras \cite{FZ} and the ring theory of the algebra of 
quantum matrices \cite{BG}.

The square case is closely related to the theory of Poisson Lie 
groups. Denote $M_n=M_{n,n}$ and $\pi_n=\pi_{n,n}$. The standard
Poisson structure on $\GLnC$, defined by Drinfeld, 
is given by
\[
\pi_{\GLn}= \sum_{i<j} 
\left( {\mathcal L}_{E_{ij}} \wedge {\mathcal L}_{E_{ji}} - 
{\mathcal R}_{E_{ij}} \wedge {\mathcal R}_{E_{ji}} \right). 
\]
Here ${\mathcal L}$ and ${\mathcal R}$ denote the left and right 
invariant vector fields on $\GLnC$, and $E_{ij}$ denote the 
standard elementary matrices considered as elements of 
$\glnC$. The canonical embedding
\[
(\GLnC, \pi_{\GLn}) \hra (M_n, \pi_n)
\]
is Poisson.

A canonical compactification of $(M_{m,n}, \pi_{m,n})$ as a 
Poisson manifold was constructed in \cite{GSV1,BGY}. Define the 
parabolic subgroup $\GL_{m+n}(\Cset)$
\[
P_{n,m}= \{ \left( \begin{smallmatrix} A & B \\ 0 & C 
\end{smallmatrix} \right) \in \GL_{m+n}(\Cset)
\mid 
A \in M_{n, n}, B \in M_{n, m}, 
C \in M_{m \times m} \}.
\]
Then $\GL_{m+n}(\Cset)/P_{n,m}$ is diffeomorphic 
to the Grassmannian $\Gr(n, m+n)$ and $M_{m,n}$ is 
embedded in it as the open Schubert cell
by
\begin{equation}
\label{emd}
X \in M_{m,n} \mt \left( \begin{smallmatrix} 
I_n & 0 \\ X & I_m 
\end{smallmatrix} \right) . P_{n,m}
\end{equation}
where $I_n$ denotes the identity matrix of size $n \times n$.
One can define the canonical Poisson structure
on $\Gr(n, m+n)$ 
\[
\pi_{\Gr(n, m+n)} = - \sum_{1\leq i < j \leq m+n} 
\chi(E_{ij}) \wedge \chi(E_{ji})
\]
where 
$\chi \colon \gl_{m+n}(\Cset) \to \Vect(\Gr(m, m+n))$ refers to 
the induced infinitesimal action from the left action
of $\GL_{m+n}(\Cset)$ on $\Gr(n, m+n)$. For a smooth
manifold $M$, $\Vect(M)$ stands for the 
Lie algebra of smooth vector fields on $M$.
The Poisson structure 
$\pi_{\Gr(n, m+n)}$ can be also defined as the
pushforward of $\pi_{\GL_{m+n}}$ 
under the natural projection 
$\GL_{m+n}(\Cset) \to \GL_{m+n}(\Cset)/P_{n,m} \cong \Gr(n, m+n)$.
The map \eqref{emd} is not Poisson but a minor correction
turns it into a Poisson map. Let us identify 
the symmetric group $S_n$ with the subgroup of 
$\GL_n(\Cset)$ consisting of permutation matrices.
Denote the longest element of $S_n$ 
by $w^n_\ci$. Then by
\cite[(3.11)]{GSV1} and \cite[Proposition 3.3]{BGY} the map
\[
i_{m,n} \colon (M_{m,n}, \pi_{m,n}) \hookrightarrow
(\GL_{m+n}(\Cset)/P_{n,m}, \pi_{GL_{m+n}}), \quad
i_{m,n}(X) = 
\left( \begin{smallmatrix} I_n & 0 \\ 
X w^n_\ci & I_m 
\end{smallmatrix} \right) . P_{n,m}
\]
is Poisson, where $I_m$ and $I_n$ are the identity matrices
of sizes $m \times m$ and $n \times n$, respectively.

Given two subsets $I \subset \{ 1, \ldots, m\}$
and $J \subset \{ 1, \ldots, n\}$ with the same 
number of elements, define by $\Delta_{I, J}(X)$ the 
minor of $X \in M_{m,n}$ (i.e. the determinant
of the submatrix of $X$) corresponding to rows
$I$ and $J$. Our first result is that 
the Hamiltonian flows of $\De_{I,J}$ are complete;
that is the corresponding integrable curves are 
defined over the full $\Cset$. In addition, we show 
that the corresponding dynamics is given by
\[
x_{ij}(t) = p(t) + \sum_{a=1}^N p_k(t) e^{ \alpha_a t}, \; 
t \in \Cset, \; \mbox{for some polynomials} \;
p(t), p_1(t), \ldots, p_N(t),
\]
such that $\deg p(t) \leq 2 |I| -1$, $\deg p_a(t) \leq 2 |I|-2$, 
$a = 1, \ldots, N$, and some $\alpha_1, \ldots, \alpha_N \in \Cset^*$.
Here $|I|$ denotes the number of elements of $I$.
The integer $N$, the polynomials $p(t)$, $p_a(t)$ and 
the exponents $\alpha_a$ depend, in general, on the 
initial conditions.
This result is proved in Section \ref{2}.

Kogan and Zelevinsky \cite{KZ} constructed 
completely integrable systems on symplectic leaves 
of the double Bruhat cells (corresponding to equal 
elements of the Weyl group) of any complex
reductive Lie group, equipped with the standard Poisson structure.
Their Hamiltonians are either twisted generalized minors 
or generalized minors. If we restrict ourselves to 
$\GLnC$, this construction produces families of 
$(n-1)(n-2)/2$ functionally independent Hamiltonians
on $M_n$ (given by minors) parametrized by 
reduced expressions of the longest element $w^n_\ci \in S_n$. 
Our result implies that, in the $A_{n-1}$ case, 
the flows of all Kogan-Zelevinsky integrable
systems are complete, even when extended to $M_n$.
In particular each of these systems provides a 
Hamiltonian action of $\Cset^{(n-1)(n-2)/2}$
on $(M_n, \pi_n)$ (and on $(\GLnC, \pi_{\GL_n})$ and 
on $({\mathrm{SL}}_n(\Cset), \pi_{\GL_n}|_{{\mathrm{SL}}_n})$). 
{\em{The action 
is analytic but not algebraic.}} This 
is carried out in Section \ref{3}. One should note 
here that, although the left and right actions of 
the standard maximal torus of
$\GLnC$ preserve $\pi_{\GLn}$,
they are not Hamiltonian. In addition, they provide 
a lot smaller abelian symmetry groups, compared to the 
constructed ones of dimension $(n-1)(n-2)/2$.

Kostant and Wallach have recently shown that all flows of the 
complexified Gelfand--Zeitlin integrable systems are complete
everywhere on $\gl_n^*(\Cset)$. (This is not true for the 
flows of the classical system on ${\mathfrak{u}}_n$.)
This paper grew out of our attempt to understand whether 
a similar fact is true for a Gelfand--Zeitlin integrable 
system for the quadratic Poisson structure $\pi_n$ 
on $M_n$. We construct a Gelfand--Zeitlin type integrable
system on $(M_n, \pi_n)$ as follows. Any projection 
from $(M_n, \pi_n)$ to $(M_k, \pi_k)$ sending 
an $n \times n$ matrix to (any) given $k \times k$ submatrix
is Poisson for all $k <n$. We get the chain of Poisson maps
\[
(M_{n}, \pi_n) \to (M_{n-1}, \pi_{n-1}) 
\to 
\ldots \to
(M_1, \pi_1)
\]
by sending a $k \times k$ matrix to its principal $(k-1) \times (k-1)$
minor for all $k = n, n-1, \ldots, 2$.
The field of rational functions on any $M_k$ is a Poisson field.
We determine its Poisson center $Z(\Cset(M_k))$ (i.e. its center 
with respect to the Poisson bracket) using a result of Kogan and Zelevinsky 
\cite{KZ}, see \prref{cent} for details. Pulling back all Poisson centers 
to $\Cset(M_n)$ produces the Poisson commutative subfield of $\Cset(M_n)$
which is the pure transcendental extension of $\Cset$ by
\begin{equation}
\label{GZham0}
\Delta_{l;k}/\Delta_{l;l-k}, 1 \leq k < l \leq n, \Delta_l,
1 \leq l \leq n.
\end{equation}
We call the set of Hamiltonians \eqref{GZham0} the Gelfand--Zeitlin  
integrable system for the quadratic Poisson bracket $\pi_n$.
Those Hamiltonians indeed define an integrable system 
on any symplectic leaf of the open double Bruhat cell in
$\GL_n(\Cset)$. Every such leaf has dimension $n(n-1)/2$.

Finally we prove that the flow of each of the Hamiltonians
\eqref{GZham0} is complete in the following sense (see \thref{complGZ}).
Denote by $D_n$ the subvariety of $M_n$ which is the union of 
zero level sets of the denominators of the Hamiltonians in \eqref{GZham0}.
Fix $h$ to be one of the Hamiltonians \eqref{GZham0} and 
$X_0 \in M_n \backslash D_n$. We prove that

{\em{there exists a curve $\gamma \colon \Cset \to M_n$ whose
entry functions are given by
\[
\sum_{p=1}^N (c^p_{ij} + d^p_{ij} t) e^{ \alpha^p_{ij} t} +
(b_{ij} + c_{ij}t + d_{ij} t^2)
\]
for some integer $N$ and 
$b_{ij}, c^{ij}, d^{ij}, c_{ij}^p, d_{ij}^p \in \Cset$, 
$\alpha_{ij}^p \in \Cset^*$
depending on $X_0$ with the following properties:

1. $\gamma^{-1}(D_n)$ is a discrete subset of $\Cset$ and

2. $\gamma \colon \Cset \backslash \gamma^{-1}(D_n) \to M_n$ is the 
integral curve of the Hamiltonian $h$ such that 
$\gamma(0) = X_0$.}}
\\ \hfill \\
{\bf Acknowledgements.} M.Y. would like to thank
the organizers of the conference Poisson 2008 for 
the invitation to participate at this very stimulating 
conference. The research of M. G. was partially supported
by  NSF grant DMS-0400484. The research of M.Y. was partially supported
by NSF grant DMS-0406057 and an Alfred P. Sloan Research 
Fellowship.
\sectionnew{Completeness of determinantal Hamiltonian flows}\lb{2}
Let
\begin{equation}
\label{IJ}
I= \{i_1 < \ldots < i_r\} \subset \{1, \ldots, m\}, \quad
J= \{j_1 < \ldots < j_r\} \subset \{1, \ldots, n\}.
\end{equation}
For $k \in \{1, \dots, n\}$ set
\[
I(i_p \to k) = (I \backslash \{i_p\}) \cup \{k\}.
\]
Recall the Poisson brackets \cite{GSV1} 
\begin{multline}
\{x_{kl}, \De_{I,J} \} = \sum_{q=1}^r \sign(i_q-k) x_{i_q l} 
\De_{I(i_q \to k),J} \\
+ \sum_{q=1}^r \sign(j_q-l) x_{k j_q} \De_{I,J(j_q\to l)},
\label{minor_br}
\end{multline}
for all $k = 1, \ldots, m$ and $l = 1, \ldots, n$
(cf. the proof of \cite[Lemma 3.2]{GSV1}).
Following \cite{GSV1}, we say that $\sign(I-k)$ is defined if 
$k \in I$, $k < i_1$, or $k > i_r$ in which cases it is 
equal to $0$, $1$, or $-1$, respectively. We recall 
the following lemma for the convenience of the reader.
\ble{commut} \cite[Lemma 3.2]{GSV1} If, in the above notation, 
$\sign(k-I)$ and $\sign(J-l)$ are defined and
\[
|\sign(I-k) + \sign(J-l)| \leq 1,
\]
then 
\[
\{x_{kl}, \De_{I,J} \} = (\sign(I-k) + \sign(J-l)) 
x_{kl} \De_{I,J}.
\]
\ele
The main result of this section is:
\bth{main} The Hamiltonian flow of any minor $\De_{I,J}$ on 
$(M_{m,n}, \pi_{m,n})$ is complete over $\Cset$. In addition, 

each matrix entry evolves quasi-exponentially according to this
flow: 
\begin{equation}
\label{evolve}
x_{kl}(t) = p(t) + \sum_{a=1}^N p_a(t) e^{ \alpha_a t}, \;
t \in \Cset, \; \mbox{for some polynomials} \;
p(t), p_1(t), \ldots, p_N(t),
\end{equation}
such that $\deg p(t) \leq 2 r-1$, $\deg p_a(t) \leq 2 r-2$, 
$a = 1, \ldots, N$, 
and some $\al_1, \ldots, \al_N \in \Cset^*$.
The integer $N$, the polynomials $p(t)$, $p_a(t)$, and 
the exponents $\alpha_a$ depend on the initial condition.
\eth
\begin{proof}
Assume that $I$ and $J$ are given by \eqref{IJ}. We will 
consider several cases for $k$ and $l$.

Case 1: $(k > i_r, l < j_1)$ or $(k< i_1, l > j_r)$ 
or $(k \in I, l \in J)$. In this
case  $\sign(I-k)$ and $\sign(J-l)$ are defined,
and $\sign(I-k) + \sign(J-l) = 0.$ According to 
\leref{commut}
\[
x_{kl}'(t) = \{x_{kl}, \De_{I,J}\}(X(t))=0
\]
and $x_{kl}(t) =x_{kl}(0)$ for all $t \in \Cset$.

Case 2: $(k \notin I, l \in J)$. The case
$(k \in I, l \notin J)$ is treated analogously.
In this case \eqref{minor_br} implies
\[
x_{kl}'(t) = \sum_{q=1}^r \sign(i_q-k) x_{i_q l}(t) 
\De_{I(i_q \to k),J}(X(t)).
\]
Taking into account case 1, we obtain that
\[
x_{kl}'(t) = c_{kl i_1} x_{k i_1}(t) + \cdots +
c_{kl i_r} x_{k i_r}(t)
\]
for some constants $c_{k i_a j_b}$ ($k \notin I$, 
$a,b =1, \ldots, r$) which depend on the 
initial condition. Thus the vectors
\[
\vec{x}_k(t) = (x_{k j_1}(t), \ldots, x_{k j_r}(t))^T, 
\quad k \notin I
\]
satisfy
\[
\vec{x}'_k(t) = C_k \vec{x}_k(t)
\]
for the constant matrices $C_k = (c_{k i_a j_b})_{a,b=1}^r$.
Therefore, in this case, the entries $x_{kl}(t)$ evolve 
according to \eqref{evolve} with 
$\deg p(t) + \deg p_1 (t) + \cdots + \deg p_N(t) \leq r-N$.

Case 3: $(k \notin I, l \notin J)$. Eq. \eqref{minor_br}
and cases 1-2 imply
\[
x_{kl}'(t) = p_0(t) +
\sum_{a=1}^N p_a(t) e^{ \alpha_a t}, \quad
t \in \Cset, 
\]
for some polynomials $p_0(t), \ldots, p_N(t)$
such that $\deg p_a(t) \leq 2 r-2$, $a=0, \ldots, N$
and some $\alpha_1, \ldots, \alpha_N \in \Cset^*$.
The statement of the theorem now follows from this.
Note that some of the cases in 1 are included in case 3, 
but the dynamics in 1 is simpler. 
\end{proof}
\sectionnew{Completeness of the Kogan--Zelevinsky integrable systems on 
$(M_n, \pi_n)$ and Hamiltonian group actions}\lb{3}
In \cite{KZ} Kogan and Zelevinsky constructed completely integrable systems 
on certain symplectic leaves of a complex simple group $G$ equipped with 
the standard Poisson structure (the leaves in all double Bruhat cells
of the type $G^{u,u}$ for some elements $u$ of the corresponding Weyl group).
We restrict our attention to $\GL_n$ and the leaves
in the open double Bruhat cell 
$\GL_n^{w^n_\ci, w^n_\ci}= B_+ w^n_\ci B_+ \cap B_- w^n_\ci B_-$
(recall that $w^n_\ci$ denotes the longest element of $S_n$).
{\em{Since $B_+ w^n_\ci B_+ \cap B_- w^n_\ci B_-$ is Zariski open 
in $M_n$ and the corresponding Hamiltonians are regular functions on
$M_n$, the Kogan--Zelevinsky integrable systems for this case
provide sets of $n(n-1)/2$ functionally independent commuting 
Hamiltonians on $(M_n, \pi_n)$.}}

Denote by $s_1, \ldots, s_{n-1}$ the simple reflections in $S_n$.
For two integers $k \leq l$ set 
\begin{equation}
\label{int}
[k,l] = \{k, k+1, \ldots, l \}.
\end{equation}
Fix a reduced expression for the longest element of $S_n$:
\begin{equation}
\label{l_expr}
w^n_\ci = s_{j_1} \ldots s_{j_{n(n-1)/2}}.
\end{equation}
Following \cite{KZ}, for $k = 1, \ldots, n(n-1)/2$ set
\[
v_k = s_{j_1} \ldots s_{j_{k-1}}, \quad u_k = (w^n_\ci)^{-1} v_k=
s_{j_{n(n-1)/2}} \ldots s_{j_k}
\] 
and define the Hamiltonians
\begin{equation}
\label{KZ}
H_k = \De_{u_k[1,k], v_k[1,k]} \in \Cset[M_n], 
\quad k = 1, \ldots, n(n-1)/2.
\end{equation}
\bth{commut} {\em{(}}Kogan--Zelevinsky, \cite{KZ}{\em{)}} 
For any reduced expression \eqref{l_expr} of the longest 
element of $S_n$, the Hamiltonians $H_1, \ldots, H_{n(n-1)/2}$ 
Poisson commute with respect to $\pi_n$.
\eth 

\thref{main} has the following immediate corollary.

\bco{Ham_act} For any reduced expression \eqref{l_expr}
of the longest element of $S_n$ the flows of the 
Kogan--Zelevinsky Hamiltonians \eqref{KZ} are complete 
on $(M_n, \pi_n)$. As a consequence from each such reduced expression
one obtains an analytic Hamiltonian action of $\Cset^{n(n-1)/2}$ 
on $(M_n, \pi_n)$.
\eco 

As was pointed out in the introduction, the left and right 
actions of the standard maximal torus of $\GLnC$ are not Hamiltonian
and provide a lot smaller abelian symmetry groups of $(M_n, \pi_n)$.
\sectionnew{Gelfand-Zeitlin integrable systems on $(M_n, \pi_n)$}\lb{4}
The field of rational functions $\Cset(P)$ of an irreducible affine 
complex Poisson variety $(P, \pi)$ is naturally a Poisson field.
Its center with respect to the Poisson bracket will be denoted 
by $Z(\Cset(P))$ and will be called the Poisson center of $\Cset(P)$.
\subsection{Poisson maps between different $(M_{m,n}, \pi_{m,n})$}
\lb{maps} 
For two subsets $I \subset \{1, \ldots, m\}$ and 
$J \subset \{1, \ldots, n\}$
we have the natural projection and embedding maps:
\[
\phi^{m,n}_{I,J} \colon M_{m,n} \to M_{|I|, |J|}, \quad
i_{I,J}^{m,n} \colon M_{|I|,|J|} \hra M_{m,n}
\]
where $|I|$ and $|J|$ denote the number of elements of $I$ and $J$.
Given $X \in M_{m,n}$, $\phi^{m,n}_{I, J}(X)$ denotes 
the submatrix of $X$ corresponding 
to rows in $I$ and columns in $J$. For $X \in M_{|I|,|J|}$, 
$i^{m,n}_{I, J}(X) \in M_{m,n}$ is the matrix whose 
$I,J$'th submatrix is $X$ and all other entries are equal to 0.

One immediately checks that:

\ble{Pmaps} For all subsets $I \subset \{1, \ldots, m\}$, 
$J \subset \{1, \ldots, n\}$ the maps
\[
\phi^{m,n}_{I,J} \colon M_{m,n} \to M_{|I|, |J|} \; \; 
{\mbox{and}} \; \;  
i^{m,n}_{I,J} \colon M_{|I|,|J|} \to M_{m,n}
\] 
are morphisms of affine Poisson spaces.
\ele
Note that if one identifies $M_n$ with $\gl^*_n(\Cset)$, then for the 
Kirillov--Kostant Poisson structure the morphisms $\phi^{n,n}_{I,J}$ are 
Poisson only if $I=J$ and $i^{n,n}_{I,J}$ are never Poisson.
\subsection{The Poisson center of $\Cset(M_n)$}
\lb{centers}
We will construct Poisson commutative subfields of $\Cset(M_n)$
by pulling back Poisson centers under a chain of projection maps.
For this purpose, we will need a description of the Poisson centers
$Z(\Cset(M_n))$. 

For $1 \leq k<l \leq n$ define $\Delta_{l;k}, \Delta'_{l;k} \in \Cset[M_n]$ by
\[
\Delta_{l;k}(X) = \Delta_{[1, k],[l-k+1, l]}(X) 
\quad \mbox{and} \quad
\Delta'_{l;k}(X) = \Delta_{[l-k+1,l],[1,k]}(X), \; \; X \in M_n.
\]
In other words $\Delta_{l;k}(X)$ and $\Delta'_{l;k}(X)$ are 
the upper--right $k \times k$ minor and the lower--left $k \times k$ 
minors of $\phi^{n,n}_{[1,l],[1,l]}(X).$

For $1 \leq l \leq n$, define also $\Delta_l \in \Cset[M_n]$
by
\[
\Delta_l(X) = \Delta_{[1,l], [1,l] }(X)= 
\det(\phi_{[1,l],[1,l]}(X)), \; \; X \in M_n.
\] 

The following proposition describes the Poisson 
centers $Z(\Cset(M_n)$.

\bpr{cent}
The (Poisson) center of the Poisson field $\Cset(M_n)$
is
\[
\Cset( \Delta_{n;1}/\Delta'_{n;n-1}, \ldots, 
\Delta_{n;n-1}/\Delta'_{n;1}, \Delta_n).
\]
\epr

It will be obtained as an easy consequence of the following 
theorem of Kogan and Zelevinsky (which is a special
case of their result \cite[Theorem 2.3]{KZ} describing 
the symplectic leaves of the standard Poisson structure
on any complex simple Lie group,
see \cite[Example 2.10]{KZ}).
\bth{thmKZ} {\em{(}}Kogan--Zelevinsky{\em{)}} The symplectic 
leaves of the Zariski open subset $\GL^{w^m_\ci, w^n_\ci}_n$
of $M_n$ are
\begin{multline}
\{ X \in \GL_n^{w^m_\ci, w^n_\ci} \mid 
\Delta_{n;1}(X)/\Delta'_{n;n-1}(X) = \const, \ldots, \\
\Delta_{n;n-1}(X)/\Delta'_{n;1}(X) = \const, 
\Delta_n(X) = \const \}.
\end{multline}
\eth  

It is clear that \thref{thmKZ} implies the inclusion
\[
Z(M_n) \supseteq \Cset( \Delta_{n;1}/\Delta'_{n;n-1}, \ldots,
\Delta_{n;n-1}/\Delta'_{n;1}, \Delta_n).
\]
The major point of \prref{cent} is the equality of the above 
fields. \\
\medskip
\noindent
{\em{Proof of \prref{cent} }}. Clearly 
$Z(\Cset(M_n))= Z(\Cset(\GL^{w^n_\ci, w^n_\ci}_n))$.
Consider the map 
\[
\psi \colon \GL^{w^n_\ci, w^n_\ci} \to (\Cset^*)^{\times n}, 
\quad \psi= (\Delta_{n;1}/\Delta'_{n;n-1}, \ldots,
\Delta_{n;n-1}/\Delta'_{n;1}, \Delta_n).
\]
Any $f \in Z(\Cset(\GL^{w^m_\ci, w^n_\ci}))$ has to be 
constant along generic symplectic leaves of 
$(\GL^{w^m_\ci, w^n_\ci}, \pi_n)$. Because of the 
Kogan--Zelevinsky theorem
such an $f$ has to be constant along the generic fibers 
of $\psi$. Now the theorem follows from the 
fact that 
\begin{equation}
\label{field}
\Cset(M_n) = 
\Cset( \Delta_{n;k}/\Delta'_{n;n-k}, 1 \leq k< n, \Delta_n, 
x_{ij}, 1< i \leq n, 1 \leq j \leq n)
\end{equation}    
(in particular $\Cset(M_n)/\psi^*\Cset( (\Cset^*)^{\times n} )$ is 
a purely transcendental extension). This is proved by 
checking that
\[
\Delta_{n;k}/\Delta'_{n;n-k} = \sum_{d=1}^k Q_{kd} x_{1,n+1-d}, \;  
k = 1, \ldots, n-1, \quad 
\Delta_n = \sum_{d=1}^n Q_{nd} x_{1,n+1-d}
\]
for some
\[
Q_{kd} \in \Cset(x_{ij}, 1<i \leq n, 1 \leq j \leq n), 1 \leq d \leq k \leq n, \quad 
Q_{kk} \neq 0.
\]
\qed
\subsection{A Gelfand--Zeitlin type integrable system}
\lb{gc}
For $k= 2, \ldots, n$ set
\[
\phi_k = \phi^{k,k}_{[1, k-1], [1, k-1]} 
\colon (M_k, \pi_k) \to (M_{k-1}, \pi_{k-1}). 
\]
We have the chain of Poisson maps
\begin{equation}
\label{chain}
(M_{n}, \pi_n) \stackrel{\phi_n}\longrightarrow (M_{n-1}, \pi_{n-1}) 
\stackrel{\phi_{n-1}}\longrightarrow 
\ldots \stackrel{\phi_2}\longrightarrow
(M_1, \pi_1),
\end{equation}
cf. \leref{Pmaps}. Following the Gelfand--Zeitlin procedure 
we induce a Poisson commutative 
subfield of $\Cset(M_n)$ by pulling back the Poisson centers 
of all fields $\Cset(M_k)$ and considering the subfield 
of $\Cset(M_n)$ generated by them. More precisely, set
\[
\phi_{n;k} = \phi_n \cdots \phi_{k+1} \colon (M_n, \pi_n) \to
(M_k, \pi_k)
\]
and consider the Poisson commutative subfield of $\Cset(M_n)$
generated by $\phi^*_{n;k}(\Cset(M_k))$ for all $k=1, \ldots n$.
Because of \prref{cent} this is the field
\[
\Cset( \Delta_{l;k}/\Delta'_{l;l-k}, 1 \leq k < l \leq n, \Delta_l, 
1 \leq l \leq n).
\] 
This proves the Poisson commutativity of the set of rational 
Hamiltonians
\begin{equation}
\label{GZham}
\Delta_{l;k}/\Delta'_{l;l-k}, 1 \leq k < l \leq n, \Delta_l,
1 \leq l \leq n
\end{equation}
on $(M_n, \pi_n)$.
We call this set the Gelfand--Zeitlin integrable system
on the matrix affine Poisson space $(M_n, \pi_n)$.
Let us note that the commutativity of the Hamiltonians 
can be also deduced from \cite[Proposition 3.4]{GSV1}. 

The set of Hamiltonians \eqref{GZham} induces 
integrable systems on the leaves of the Zariski open
subset $\GL_n^{w^n_\ci, w^n_\ci}$ as follows. Among those Hamiltonians, 
$n$ are Casimirs for $(\GL_n^{w^n_\ci, w^n_\ci}, \pi_n)$ and restrict 
to constants, cf. \thref{thmKZ}. The remaining $n(n-1)/2$ are 
functionally independent on every symplectic leaf of 
$\GL_n^{w^n_\ci, w^n_\ci}$ because iterating \eqref{field}
gives
\[
\Cset(M_n) =
\Cset( \Delta_{l;k}/\Delta'_{l;l-k}, 1 \leq k< l \leq n, 
\Delta_l, 1 \leq l \leq n,
x_{ij}, 1 \leq j < i \leq n).
\]
Since all leaves of $\GL_n^{w^n_\ci, w^n_\ci}$ have dimension
$n(n-1)$ this induces integrable systems on them.   

The integrable systems \eqref{GZham} are related to the Gelfand--Zeitlin
integrable systems defined by Flaschka and Ratiu \cite{FR, FR2} in the
sense that they are produced by similar procedures but live on dual 
Poisson Lie groups. The one of Flaschka and Ratiu is on the 
dual Poisson--Lie group $U_n^*$ of $U_n$ (and thus its complexification 
is a Hamiltonian system on the dual Poisson Lie group of 
$\GL_n(\Cset)$). The Hamiltonians \eqref{GZham} define dynamics on 
$M_n$ and on its Poisson submanifold $\GL_n(\Cset)$.
It was 
conjectured by Flaschka and Ratiu \cite{FR, FR2}, and proved by Alekseev and 
Meinrenken \cite{AM} that the Gelfand--Zeitlin integrable systems 
on $U^*_n$ are globally isomorphic to the classical ones. 
At the same time the systems \eqref{GZham} are not isomorphic to 
Hamiltonian systems
for linear Poisson structures since the standard Poisson structure
on $\GL_n(\Cset)$ cannot be globally linearized. 

\bre{different} Instead of the chain of maps \eqref{chain}
one can consider any chain of projections $\phi^{k,k}_{I,J}$
where $|I|=|J|$. Clearly there are $n!^2$ different chains 
of this type. For each of them one can obtain a set of $n(n-1)/2$
functionally independent commuting Hamiltonians on $(M_n, \pi_n)$.
The resulting integrable systems are different since the right 
and left actions of $S_n$ on $M_n$ (permuting rows and columns) 
do not preserve the Poisson structure $\pi_n$.  
\ere
\subsection{Completeness of the flows of the  Gelfand--Zeitlin integrable 
system for $(M_n, \pi_n)$}
\lb{complgz}
Denote the singular locus of the Hamiltonians \eqref{GZham} 
\[
D_n = \{ X \in M_n \mid
\Delta'_{l; k } (X) = 0 \; \; \mbox{for some} \; \; 
1 \leq k < l \leq n\}.
\] 

We prove that the flows of the Hamiltonians \eqref{GZham} are complete 
in the following sense.

\bth{complGZ} 
Let $h$ be one of the Hamiltonians \eqref{GZham} and 
$X_0 \in M_n \backslash D_n$. 

There exists a curve $\gamma \colon \Cset \to M_n$ whose
entry functions are given by
\begin{equation}
\label{coeff}
\sum_{p=1}^N (c^p_{ij} + d^p_{ij} t) e^{ \alpha^p_{ij} t} +
(b_{ij} + c_{ij}t + d_{ij} t^2)
\end{equation}
for some integer $N$ and 
$b_{ij}, c^{ij}, d^{ij}, c_{ij}^p, d_{ij}^p  \in \Cset$,
$\alpha_{ij}^p \in \Cset^*$
depending on $X_0$ with the following properties:

1. $\gamma^{-1}(D_n)$ is a discrete subset of $\Cset$ and

2. $\gamma \colon \Cset \backslash \gamma^{-1}(D_n) \to M_n$ is the 
integral curve of the Hamiltonian $h$ such that 
$\gamma(0) = X_0$.
\eth  
\begin{proof}
We will prove the statement for the case 
$h = \Delta_{l;k}/\Delta'_{l;l-k}$ for some 
$1 \leq k < l \leq n$. The case of $h = \Delta_l$ 
is simpler and is left to the reader.
 
1. Let $i,j \in [1,l]$. Then 
\[
\{ x_{ij}, \Delta_{l;k}/\Delta'_{l;l-k} \} = 0
\]
because the map $\phi^{n,n}_{[1,k],[1,k]}$ is Poisson and we 
define
\begin{equation}
\label{e1}
\gamma_{ij}(t) = (X_0)_{ij} = \const
\end{equation}
for those $i,j$'s.

2. Let $i \in [1,k]$, $j \in [l+1, n]$
or $i \in [l+1,n]$, $j \in [1,l-k]$. We set 
$\epsilon= - 1$ and $1$, respectively. 
From \leref{commut} one gets
\[
\{ x_{ij}, \Delta_{l;k}/\Delta'_{l;l-k} \} = \epsilon\ x_{ij}
\Delta_{l;k}/\Delta'_{l;l-k}.
\]
Define 
\begin{equation}
\label{e2}
\gamma_{ij}(t) = (X_0)_{ij} 
e^{\epsilon \Delta_{l;k}(X_0)/\Delta'_{l;l-k}(X_0)}.
\end{equation}
For those values of $i$ and $j$,
this is the evolution of the matrix entries of 
the integral curve of $h$ through $X_0$ outside the singular 
locus $D_n$.

3. Let $i \in [k+1, l]$, $j \in [l+1, n]$. 
Using \eqref{minor_br} and \leref{commut} we get
\begin{multline}
\label{mult1}
\{ x_{ij}, \Delta_{l;k}/\Delta'_{l;l-k} \} = 
- \Big( \sum_{q=1}^k x_{q j}\Delta_{[1,k](q \to i), [l-k+1, l]} 
+
\\ 
\sum_{q=l-k+1}^l x_{iq} \Delta_{[1,k], [l-k+1, l](q\to j)} 
- \Delta_{l;k} \Big)
/ \Delta'_{l; l-k}.
\end{multline}
Substituting $x_{ab}$ with $\gamma_{ab}(t)$ from \eqref{e1}--\eqref{e2} 
in the RHS of \eqref{mult1}, gives that outside
$D_n$ the matrix entries of the integral curve of $h$ through $X_0$ 
satisfy
\begin{equation}
\label{i3}
x'_{ij}(t) = \sum_{p=1}^N c^p_{ij} e ^{ \alpha^p_{ij} t} + d_{ij}
\end{equation}
and thus
\begin{equation}
\label{e3}
x_{ij}(t) = \sum_{p=1}^N c^p_{ij} e ^{ \alpha^p_{ij} t} +
(c_{ij} + d_{ij} t)
\end{equation}
for some integer $N$ and $c^{ij}, d^{ij}, c_{ij}^p \in \Cset$,
$\alpha_{ij}^p \in \Cset^*$.
For $i \in [k+1, l]$, $j \in [l+1, n]$ define 
$\gamma_{ij}(t)$ by the RHS of \eqref{e3}.

4. Let $i \in [l+1, n],$ $j \in [l-k+1,l]$. In this case one computes
the Poisson bracket $\{ x_{ij}, \Delta_{l;k}/\Delta'_{l;l-k} \}$ 
similarly to \eqref{mult1} and proves that it is given by a 
formula of the type \eqref{i3} where again 
$c^{ij}, d^{ij}, c_{ij}^p \in \Cset$, $\alpha_{ij}^p \in \Cset^*$ 
depend on $X_0$.
For $i \in [l+1, n],$ $j \in [l-k+1,l]$ we define $\gamma_{ij}(t)$ by the 
RHS of the latter formula. That formula also gives the evolution of 
the matrix entries of the integral curve of $h$ through 
$X_0$ outside the singular locus $D_n$.

5. Let $i,j \in [l+1,n]$. Using \eqref{minor_br} one obtains
\begin{multline}
\label{mult2}
\{ x_{ij}, \Delta_{l;k}/\Delta'_{l;l-k} \} = 
- \Big( \sum_{q=1}^k x_{q j}\Delta_{[1,k](q \to i), [l-k+1, l]} 
+
\\ 
\sum_{q=l-k+1}^l x_{iq} \Delta_{[1,k], [l-k+1, l](q\to j)} \Big)
/ \Delta'_{l; l-k}+ 
\Delta_{l;k} \Big( \sum_{q=1}^k x_{q j}\Delta_{[k+1,l](q \to i), [1,l-k]} 
+
\\ 
\sum_{q=l-k+1}^l x_{iq} \Delta_{[k+1,l], [1, l-k](q\to j)} \Big)
/ (\Delta'_{l; l-k})^2.
\end{multline}
We substitute $x_{ab}$ with $\gamma_{ab}(t)$ from the formulas in cases 
1--4 in the RHS of \eqref{mult2}, and get that outside
$D_n$ the matrix entries of the integral curve of 
$h$ through $X_0$ satisfy
\begin{equation}
\label{i5}
x'_{ij}(t) = \sum_{p=0}^N (c^p_{ij} + d^p_{ij} t) e ^{ a^p_{ij} t}
\end{equation}
for some $c_{ij}^p, \alpha_{ij}^p \in \Cset$.
We will assume that $\alpha_{ij}^0=0$ and set $c^0_{ij}= d^0_{ij}=0$ if 
polynomial terms are not present. Thus, outside $D_n$ the 
integral curve of $h$ through $X_0$ is integrated to
\begin{equation}
\label{e5}
x_{ij}(t) = \sum_{p=1}^N (c^p_{ij} + d^p_{ij} t) e^{ a^p_{ij} t} +
(b_{ij} + c_{ij}t + d_{ij} t^2)
\end{equation}
for some $b_{ij}, c_{ij}, d_{ij}, c_{ij}^p, d_{ij}^p \in \Cset$,
$\alpha_{ij}^p \in \Cset^*$. 
For $i, j \in [l+1, n]$ we define
$\gamma_{ij}(t)$ by the RHS of \eqref{e5}.

The cases 1--5 define an analytic curve $\gamma \colon \Cset \to M_n$ 
with the property \eqref{coeff}. If $H$ is any of the Hamiltonians
\eqref{GZham} and $\Delta$ is its denominator, 
then $\Delta(\gamma(t))$ is an entire function which does not 
vanish identically since $\Delta(\gamma(0))= \Delta(X_0) \neq 0$. 
Therefore 
\begin{equation}
\label{gamma-1}
\gamma^{-1} ( \{ X \in M_n \mid \Delta(X) = 0 \} ) 
\end{equation}
is a discrete subset of $\Cset$. Thus $\gamma^{-1}(D_n)$ 
is a discrete subset of $\Cset$ since it is 
the union of all sets \eqref{gamma-1} where $H$ runs over 
the Hamiltonians \eqref{GZham}. Finally, outside 
$\gamma^{-1}(D_n)$, $\gamma(t)$ is the integral curve 
of the Hamiltonian $h$ through $X_0$ because of the way
the matrix entries of $\gamma(t)$ were constructed.   
\end{proof}

\end{document}